\documentclass[11pt]{amsart}
\usepackage{epsfig}
\usepackage{amsmath}
\usepackage{amsfonts}
\usepackage{amssymb}
\usepackage{color}

\title{The image of a finely holomorphic map is pluripolar}

\author{Armen Edigarian, Said El Marzguioui and Jan Wiegerinck\\}

\baselineskip

\newcounter{bean}

\newtheorem{theorem}{Theorem}[section]
\newtheorem{corollary}[theorem]{Corollary}
\newtheorem{lemma}[theorem]{Lemma}
\newtheorem{proposition}[theorem]{Proposition}
\theoremstyle{definition}
\newtheorem{definition}[theorem]{Definition}
\newtheorem{remark}[theorem]{Remark}
\newtheorem{example}[theorem]{Example}

\newcommand{\CC}{\mathbb{C}}

\DeclareMathOperator*{\f-limsup}{f-\,limsup}
\DeclareMathOperator*{\PSH}{PSH}
\begin{document}

\maketitle \footnote{2000 Mathematics Subject Classification
32U15, 31C40} \footnote{The work is part of the Research Grant No.
1PO3A 005 28, which is supported by public means in the programme
for promoting science in Poland in the years 2005-2008.}
\begin{abstract}
We prove that the image of a finely holomorphic map on a fine
domain in $\mathbb{C}$ is a pluripolar subset of $\mathbb{C}^{n}$.
We also discuss the relationship between pluripolar hulls and
finely holomorphic functions.
\end{abstract}
\begin{scriptsize}
\end{scriptsize}

\section{Introduction}

A subset $E \subset \mathbb{C}^{n}$ is said to be
\emph{pluripolar} if for each point $a \in E$ there is an open
neighborhood $\Omega$ of $a$ and a function $\varphi $ (
$\not\equiv -\infty$) plurisubharmonic in $\Omega$,
($\varphi\in\PSH(\Omega)$) such that
$$
E\cap \Omega \subset \{ z\in \Omega :  \varphi (z)=- \infty\}.
$$
It is a fundamental result of Josefson \cite{Jo78} that this local
definition is equivalent to the global one, i.e., in this
definition one can assume $\varphi$ to be plurisubharmonic in all
of $\mathbb{C}^{n}$ with
$$
E\subset \{ z\in \mathbb{C}^{n} : \varphi (z)=- \infty\}.
$$
$E$ is called \emph{complete pluripolar} (in $\mathbb{C}^{n}$) if
for some plurisubharmonic function $\varphi \in
\text{PSH}(\mathbb{C}^{n})$, we have $E=\{z \in \mathbb{C}^{n}:
\varphi(z) = -\infty\}$. Unlike the situation in classical
potential theory, pluripolar sets often "propagate"; it may happen
that any PSH function $\varphi$ which is $-\infty$ on a pluripolar
set $E$ is automatically $-\infty$ on a larger set. For example,
if the $-\infty$ locus of a PSH function $\varphi$ contains a
non-polar piece of a complex analytic variety $A$, then the set
$\{z \in \mathbb{C}^{n}: \varphi(z) = -\infty \}$ must contain all
the points of $A$. However, the structure of pluripolar sets may be much more complicated,
cf. \cite{Le88,CLP}.
%by a suitable choice of parameters in
%Wermer's famous example (cf. \cite{We82}), Levenberg (see
%\ constructed an example of a compact non-complete
%pluripolar set which hits every complex analytic variety in a
%polar set. More recently, Coman, Levenberg, and Poletsky (see
%\cite{}) have constructed a non-pluripolar set which intersects
%every complex analytic disc in a discrete set. These two important
%results reflect the complicated nature of the structure of
%pluripolar sets and the curious phenomenon of propagation they
%exhibit. In recent years, c
Completeness of pluripolar sets has received growing attention,
and in particular cases many results were obtained, see
\cite{CLP,EW03,EW04,E,EJ,Le88,Si03,Wi00,Zw05}. But our knowledge
and understanding of the general situation is fragmentary, and a
good characterization of complete pluripolar sets is still
lacking, even in the case of the graph of an analytic function.

Recently, in \cite{EJ} Edlund and J\"oricke have connected the
propagation of the graph of a holomorphic function as a
pluripolar set to fine analytic
continuation of the function.
\begin{theorem}[Edlund and J\"oricke, \cite{EJ} Theorem 1] \label{thm1} Let $f$ be holomorphic in
the unit disc $\mathbb{D} \subset \mathbb{C}$ and let $p\in
\partial \mathbb{D}$. Suppose that $f$ has a finely holomorphic
continuation $F$ at $p$ to a closed fine neighborhood $V$ of $p$.
Then there exists another closed fine neighborhood $V_{1} \subset
V$ of $p$, such that the graph $\Gamma_{F}(V_{1})$ is contained in
the pluripolar hull of $\Gamma_{f}(\mathbb{D})$.
\end{theorem}

The definition of the pluripolar hull and necessary preliminaries
about finely holomorphic functions are presented in Section 2.

In view of this result, it is reasonable to try and investigate
the connection between finely holomorphic functions and pluripolar
sets. Using some of Fuglede's fundamental results in fine
potential theory, we can easily prove stronger results. Moreover,
our method allows to give shorter proofs of known results about
pluripolar hulls. \\

%\begin{theorem} Let $f$ be holomorphic in a connected open set $U \subset \mathbb{C}$
%and let $p\in \partial U$. Suppose that $f$ has a finely
%holomorphic continuation $F$ at $p$ to a finely open and finely
%connected neighborhood $V$ of $p$. Then $ \Gamma_{F}(V) \subset
%(\Gamma_{f}(U))^{\ast}_{\mathbb{C}^{2}}$. Moreover, if $E$ is a
%non-polar subset of $V\cap \overline{U}$ then $ \Gamma_{f}(U) \cup
%\Gamma_{F}(V)\subset (\Gamma_{f}(E))^{\ast}_{\mathbb{C}^{2}}$.
%\end{theorem}

%The next results are the main results of this paper. The proofs
%require two lemmas and will be given in Section 3.
%\begin{proposition} Let $f$ :\ $U$ $\longrightarrow$ $\mathbb{C}^{n}$, $f(z)=(f_{1}(z),..., f_{n}(z))$, be a finely
%holomorphic map on a finely open subset $U\subseteq \mathbb{C}$.
%Then the graph $\Gamma _{f}(U)$ of $f$ is a pluripolar subset of
%$\mathbb{C}^{n+1}$.
%\end{proposition}

%As a consequence of the inversion theorem for finely holomorphic
%functions (cf. (\cite{Fu81}, Theorem 13)) we can put Proposition
%1.3 in the following more general form

Our main results is the following theorem. Theorem \ref{thm1} is a special case of it.

\begin{theorem} \label{thm2} Let $f$ :\ $U$ $\longrightarrow$ $\mathbb{C}^{n}$, $f(z)=(f_{1}(z),..., f_{n}(z))$, be a
finely holomorphic map on a finely open subset $U\subseteq
\mathbb{C}$. Then the image $f(U)$ of $U$ is a pluripolar subset
of $\mathbb{C}^{n}$. Moreover, if $E$ is a non polar subset of
$U$, then the pluripolar hull of $f(E)$ contains $f(U)$.
\end{theorem}

%For example, suppose that $K \subset\mathbb{C}$ is a compact set
%with non-empty fine interior $K'$. Every function $f \in R(K)$
%(the uniform closure of the algebra of restrictions to $K$ of
%holomorphic functions in open sets containing $K$) is finely
%holomorphic in $K'$ (cf. \cite{Fu81}, page 75). Hence, by the
%above proposition the graph $\Gamma_{f}(K') =\{(z,f(z)): z\in K'
%\}$ is a pluripolar subset of $\mathbb{C}^{2}$.
Note that in general $U$ may not have any Euclidean interior points.
The theorem applies e.g. to Borel-type series like
\begin{equation}
f(z)=\sum_{j=1}^\infty\frac{c_j}{2^j(z-a_j)},
\end{equation}
where $c_j$ are very small and $\{a_j\}$ is dense in $\CC$. We
will elaborate such an example in Section \ref{secExample}.

%For $n=1$, a partial converse of Proposition 1.3 was proved by
%Tomas Edlund in his thesis \cite{E}. Namely, he proved that if $f$
%is a function of class $C^{2}$ on a finely open set $V$, and the
%graph $\Gamma_{f}(V)$ of $f$ is pluripolar, then $f$ is finely
%holomorphic in $V$. Edlund's result together with our Proposition
%1.3 give actually (in the particular case $n=1$) a partial "fine"
%analog of the important theorem of N. V. Shcherbina that was
%obtained shortly before (see \cite{Sh03}). Shcherbina's results
%asserts that the graph $\Gamma_{f}(\Omega)$ of a continuous
%function $f$ on an open set $\Omega \subset \mathbb{C}^{n}$ is
%pluripolar subset of $\mathbb{C}^{n+1}$ if and only if $f$ is
%holomorphic. It is therefore a natural question to ask whether the
%$C^{2}$-regularity in Edlund's theorem can be weakened to just
%fine continuity.

The next theorem is a simple, precise,
and complete interpretation of recent results of
the first and the third author (see \cite{EW03,EW04}).

\begin{theorem} \label{thm3} Let $D$ be a domain in $\mathbb{C}$ and let $A$ be a closed polar subset of $D$.
Suppose that $f \in \mathcal{O}(D \backslash A)$ and that
$z_{0}\in A$. Then the following
conditions are equivalent: \\
(1) $(\{z_{0}\}\times \mathbb{C})\cap (\Gamma_{f})^{\ast}_{D \times \mathbb{C}} \neq \emptyset$.\\
(2) $f$ has a finely holomorphic extension $\tilde f$ at $z_{0}$. \\
Moreover, if one of these conditions is met, then
$(\{z_{0}\}\times \mathbb{C})\cap (\Gamma_{f})^{\ast}_{D \times
\mathbb{C}}= (z_{0}, \tilde{f}(z_{0}))$.
\end{theorem}

The proofs of the above results are given in Section 3. Our
arguments rely heavily on results from fine potential theory.
Since this theory is not of a very common use in the study of
pluripolar sets, we will recall some basic facts about it. This is
done in Section 2. Using the same ideas as in the proof of Theorem
1.2 and recent results on finely plurisubharmonic functions, cf.
\cite{E-W2} we will prove in Section 5 a version of Theorem 1.1
for functions of several variables. In Section 6 we discuss some
consequences of Theorem \ref{thm2} and some open problems.\\
\textbf{Acknowledgments}. Part of this research was carried out
while the second author was visiting the mathematics department at
Copenhagen university; he would like to thank the department for
its hospitality, and express his gratitude to Professor Bent
Fuglede for his invitation and several helpful and interesting
discussions.

\section{Preliminaries}
\subsection{Pluripolar hulls}
Let $E$ be pluripolar set in $\mathbb{C}^{n}$. The
\emph{pluripolar hull} of $E$ relative to an open subset $\Omega$
of $\mathbb{C}^{n}$ is the set
$$
E_{\Omega}^{\ast}=\{z\in \Omega : \text{for} \ \text{all} \
\varphi \in \PSH (\Omega) : \varphi |_{E}= - \infty \
\Longrightarrow \varphi(z)= - \infty\}.
$$
The notion of the pluripolar hull was first introduced and studied
by Zeriahi in \cite{Ze89}. The paper \cite{LP99} of Levenberg and
Poletsky
contains a more detailed study of this notion.\\
Let $f$ be a holomorphic function in an open set $\Omega \subseteq
\mathbb{C}^{n}$. We denote by $\Gamma_{f}(\Omega)$ the graph of
$f$ over $\Omega$,
$$
\Gamma_{f}(\Omega)= \{(z, f(z)) : z\in \Omega\}.
$$
It is immediate that $\Gamma_{f}(\Omega)$ is a pluripolar subset
of $\mathbb{C}^{n+1}$. The pluripolar hull of the graph of a
holomorphic functions was studied in several papers. (See
\cite{EW03,EW04,E,EJ,Si03,Wi00,Zw05}).

Of particular interest for our present considerations is the
following (see \cite{EW03,EW04}).
\begin{theorem}[Edigarian and Wiegerinck] \label{thm4}Let $D$ be a domain in $\mathbb{C}$ and let $A$ be
a closed polar subset of $D$. Suppose that $f \in \mathcal{O}(D
\backslash A)$ and that $z_{0} \in A$. Then the following
conditions are equivalent:\\
(1) $(\{z_{0}\}\times \mathbb{C})\cap (\Gamma_{f})^{\ast}_{D \times \mathbb{C}} \neq \emptyset$.\\
(2) the set $\{z\in D \backslash A : |f(z)| \geq R \}$ is thin at
$z_{0}$ for some $R > 0$.
\end{theorem}
\subsection{Fine potential theory}
In this subsection we gather some definitions and known results
from fine potential theory that we
will need later on.\\
The {\em fine topology} on an open set $\Omega$ is the weakest
topology on $\Omega$ making all subharmonic functions continuous.
If $\Omega_1\subset \Omega_2$ are domains, then the fine topology
on $\Omega_1$ coincides with the restriction of the fine topology
in $\Omega_2$ to $\Omega_1$. The following results, except 4)
which is obvious, are due to Fuglede and can be found in
\cite{Fu72}, Chapter III-IV.
\begin{proposition}\label{prop0}
1) The fine topology is locally connected.\\
2) Every usual domain is also a fine domain. \\
3) If $U$ is a fine domain and $E$ is a polar set, then
$U\setminus E$ is a fine domain, in particular it is connected.\\
4) The fine topology has a neighborhood basis consisting of fine
neighborhoods that are Euclidean compact.
\end{proposition}

The fine topology has no infinite compact sets and is not
Lindel\"{o}f. However, the following property can serve as a
replacement. (see e.g. \cite{DO84}, page 181).
\begin{theorem}(Quasi-Lindel\"{o}f property)\label{thm5} An arbitrary union of finely open
subsets of $\mathbb{C}$ differs from a suitable countable subunion
by at most a polar set.
\end{theorem}

We now formulate the definitions and results concerning fine potential theory,
that we will use in the present paper. All of these, the proofs, and much more can be found in \cite{Fu72}.
All are quite natural in comparison with the classical situation. First we give the definitions.

\begin{definition} \label{def1}
A function $\varphi$ :\ $U$ $\longrightarrow$ $[-\infty, +\infty[$ defined on a
finely open set $U \subseteq \mathbb{C}$ is said to be finely
hypoharmonic if $\varphi$ is finely upper semicontinuous and if
$$
\varphi(z) \leq \int \varphi
d\varepsilon_{z}^{\mathbb{C}\backslash V}, \ \forall z \in V \in
\mathcal{B}(U).
$$
(It is part of the requirement that the integral exists).
$\varphi$ is finely subharmonic if, moreover, $\varphi$ is finite
on a finely dense subset of $U$.
\end{definition}

Here $\mathcal{B}(U)$ denotes the class of all finely open sets
$V$ of compact closure $\overline {V}$ (in the usual topology)
contained in $U$, and $\varepsilon_{z}^{\mathbb{C}\backslash V}$
is the swept-out of the Dirac measure $\varepsilon_{z}$ onto
${\mathbb{C}\backslash V}$. It is carried by the fine boundary
$\partial_{f}V$ of $V$. This swept-out
measure boils down to the usual harmonic measure if $V$ is a usual
open set.
\begin{theorem}\label{thm7}
% 1) Any usual subharmonic
% function is finely subharmonic on its domain.\\
1)Finely subharmonic functions on a finely open set $\Omega$ form
a convex cone that is stable  under pointwise supremum for finite
families, and closed under finely locally uniform
convergence.\\
2) A pointwise infimum of a lower directed family of finely
subharmonic functions in a fine domain $\Omega $ is either finely
subharmonic, or it identically
equals $-\infty$. \\
3) A finely subharmonic function $f$ on a finely open set $\Omega$
has a finely subharmonic restriction to every finely open subset
of $\Omega$. Conversely, suppose that $f$ is finely subharmonic in
some fine neighborhood of each point of $\Omega$. Then $f$ is
finely subharmonic in $\Omega$, i.e., finely subharmonic functions
have the {\em sheaf property}.
\end{theorem}
Next we mention the results.

\begin{proposition}\label{prop1}\cite{Fu74}
In a usual open set in $\mathbb{C}$ finely subharmonic functions
are just subharmonic ones, and the restriction of a usual
subharmonic function to a finely open set is finely subharmonic.
\end{proposition}

\begin{theorem} \label{thm6} Let $h$ :\ $U$ $\longrightarrow$ $[-\infty,+\infty[$ be a finely hypoharmonic function
on a fine domain $U \subset \mathbb{C}$. Then either the set $\{z
\in U : h(z) = - \infty\}$ is a polar subset of $U$ and $h$ is
finely subharmonic, or $h \equiv -\infty$.
\end{theorem}

\subsection{Finely holomorphic functions}

Shortly after that fine potential theory was established, several
authors turned their attention to developing the analog of
holomorphic functions on a fine domain. See \cite{Fu81},
\cite{Fu88} and the references therein. Fuglede's paper
\cite{Fu81} is our main reference for what follows.

For a compact set $K$ in $\CC$, we denote by $R(K)$ the uniform
closure on $K$ of the set of rational functions with poles outside
$K$. By Runges's theorem one can just as well take the closure of
the set of functions holomorphic in a neighborhood of $K$.

\begin{definition}\label{def2} Let $U$ be a finely open set in ${\mathbb{C}}$. A function
$f$ :\ $U$ $\longrightarrow$ $\mathbb{C}$ is called finely
holomorphic, if every point of $U$ has a compact (in the usual
topology) fine neighborhood $K \subset U$ such that the
restriction $f\mid_{K}$ belongs to $R(K)$.
\end{definition}

As we shall see below, finely holomorphic functions share many properties with ordinary holomorphic functions.
We will now assemble the results which we will need in the sequel.

\begin{theorem} \label{thm8}
A function $f$ :\ $U$ $\longrightarrow$ $\mathbb{C}$ defined in a
finely open set $U \subseteq \mathbb{C}$ is finely holomorphic if
and only if every point of $U$ has a fine neighborhood $V
\subseteq U $ in which f coincides with the Cauchy-Pompeiu
transform of some compactly supported function $ \varphi \in
L^{2}( \mathbb{C}) $ with $ \varphi =0 $ a.e. in V:
$$
f(z)= \int \limits_{\mathbb{C}} \frac{1}{z-\zeta} \varphi (\zeta)
d \lambda (\zeta), \ \ z\in V.
$$
\end{theorem}

\begin{theorem}\label{thm9}
A finely holomorphic function on a Euclidean open set is
holomorphic in the usual sense.
\end{theorem}

\begin{theorem}\label{thm10}
1) A finely holomorphic function $f$ on a fine
domain has at most countably many zeros (unless $f \equiv 0$).\\
2) A finely holomorphic function $f$ is infinitely finely
differentiable, and all its fine derivatives $f^{(n)}$ are finely
holomorphic.\\
3) Let $f$ be a finely holomorphic in a finely open set $U \subset
\mathbb{C}$. Suppose that the fine derivative $f'$ of $f$ does not
vanish at some point $z_{0} \in U$. Then one can find a finely
open neighborhood $W \subseteq U$ of $z_{0}$ such that $f|_{W}$ :\
$W$ $\longrightarrow$ $f(W)$ is bijective and the inverse function
$f^{-1}$ is finely holomorphic in the finely open set
$f(W)$.\\
4) The composition of finely holomorphic functions is finely holomorphic where it is defined.\\
5) Let $U$ be finely open and $z_0\in U$. If $f$ is finely
holomorphic on $U\setminus \{z_0\}$ and bounded in a punctured
fine neighborhood of $z_0$, then $f$ extends as a finely
holomorphic function to $U$.
\end{theorem}

\section{Pluripolarity of finely holomorphic curves}

A {\em finely holomorphic curve} is a pair $(U,f)$ where $U$ is a
fine domain and $f=(f_1,\dots,f_n) : U\to\CC^n$ is a finely
holomorphic map. As usual we will identify a curve with its image.

\begin{lemma} \label{lem1} Let $U\subseteq\mathbb{C}$ be a fine domain, and let $f$ :\
$U$ $\longrightarrow$ $\mathbb{C}^{n}$, $f(z)=(f_{1}(z)\ldots,
f_{2}(z))$, be a finely holomorphic map. Suppose that $h$ :\
$\mathbb{C}^{n}$ $\longrightarrow$ $[-\infty,\ +\infty [$ is a
plurisubharmonic function. Then the function $h \circ f$ is either
finely subharmonic on $U$ or $\equiv -\infty$.
\end{lemma}

\begin{proof}
 First, we assume that $h$ is everywhere finite and
continuous. Let $a \in U$. Definition \ref{def2} gives us
 a compact (in the usual topology) fine neighborhood
$K$ of $a$ in $U$, and $n$ sequences $(f_{j}^{k})_{k \geq 0 }$,
$j= 1,\ldots,n$, of holomorphic functions defined in Euclidean
neighborhoods of $K$ such that
$$
f_{j}^{k}|_{K} \longrightarrow f_{j}|_{K}, \ j=1,\ldots,n \
\text{uniformly}.
$$

Clearly, $(f_{1}^{k},\ldots, f_{n}^{k})$ converges uniformly on
$K$ to $(f_{1},\ldots, f_{n})$. Since $h$ is continuous, the
sequence $h(f_{1}^{k},\ldots, f_{n}^{k})$, of finite continuous
subharmonic functions, converges uniformly to $h(f_{1},\ldots,
f_{n})$ on $K$. According to Theorem \ref{thm7} 1),
$h(f_{1},\ldots,f_{n})$ is finely
subharmonic in the fine interior of $K$.\\
Suppose now that $h$ is arbitrary. We can assume that the fine
interior of $K$ is finely connected. Let $(h_{m})_{m \geq 0}$ be a
decreasing sequence of continuous plurisubharmonic functions which
converges (pointwise) to $h$. By the first part of the proof,
$h_{m}(f_{1},\ldots, f_{n})$ is a decreasing sequence of finely
subharmonic functions in the fine interior of $K$. The limit
function $h(f_{1},\ldots, f_{n})$ is by Theorem \ref{thm7} 2)
 finely subharmonic or identically $-\infty$
in the fine interior of $K$.  The sheaf property (Theorem
\ref{thm7} 3)) implies that $h(f_{1},\ldots, f_{n})$ is indeed
finely subharmonic in all of $U$ or is identically equal to
$-\infty$.
\end{proof}
\begin{remark} \label{rema1} The above lemma was also independently proved by
Fuglede.
\end{remark}
\begin{lemma} \label{lem2} Let $f$ :\ $U$ $\longrightarrow$ $\mathbb{C}^{n}$, $f(z)=(f_{1}(z),\ldots, f_{n}(z))$,
be a finely holomorphic map on a fine domain $U \subset
\mathbb{C}$ which contains a disc with positive radius. Then
$f(U)$ is a pluripolar subset of $\mathbb{C}^{n}$.
\end{lemma}

\begin{proof} Let $D(a, \delta) \subset U$ be a small disc in
$U$. Since $f$ is a holomorphic map on $D(a, \delta)$ (Theorem \ref{thm9}),
$f(D(a, \delta))$ is a pluripolar subset
of $ \mathbb{C}^{n}$. By Josefson's theorem there exists a
plurisubharmonic function $h \in\PSH(\mathbb{C}^{n})$
($\not\equiv -\infty$) such that $ h(f_{1}(z),\ldots, f_{n}(z))= -
\infty $, $\forall z \in D(a, \delta)$. According to Lemma \ref{lem1},
the function $g(z)= h(f_{1}(z),\ldots, f_{n}(z))$ is finely
subharmonic on $U$ or $\equiv -\infty$. Since it assumes $
-\infty$ on a non polar subset of $U$, it must be identically
equal to $- \infty$ on $U$ by Theorem \ref{thm6}. Hence $h|_{f(U)}=-
\infty $, and $f(U)$ is, therefore, pluripolar.
\end{proof}

%Using now Lemma 3.1 combined with the same arguments as in the
%above proof, we can give a simple proof of theorem 1.2.

%\textit{Proof of Theorem 1.2}. Denote by $g$ the finely
%holomorphic function which is equal to $f$ on $U$ and to $F$ on
%$V$. Let $h\in$ PSH$(\mathbb{C}^{2})$ be a plurisubharmonic
%function such that $h(z, f(z))= -\infty , \ \forall z \in U$.
%According to Lemma 3.1, the function $z \rightarrow h(z, g(z))$ is
%finely subharmonic on $U \cup V$ or $\equiv -\infty$. Moreover, it
%assumes $-\infty$ on the non-polar set $U$. Since $U\cup V$ is a
%fine domain, Theorem 2.4 asserts that $h(z, g(z))$ must be
%identically $- \infty$ on $U \cup V $. Hence $ \Gamma_{F}(V)
%\subset (\Gamma_{f}(U))^{\ast}_{\mathbb{C}^{2}}$. The second
%statement can be
%proved similarly. See Proposition 4.1 below for a more general results. \qed \\
%\textbf{Remark 2} The proof of Theorem 1.1 given by Edlund and
%J\"oricke in \cite{EJ} uses rather complicated harmonic measure
%estimates. In fact, the harmonic measure (especially, the two
%constant theorem) is the main ingredient in the study of
%pluripolar hulls. Its use has become quite standard. However,
%Lemma 3.1 combined with Theorem 2.4 provides an efficient
%alternative of the harmonic measure in some situations.

\begin{proposition} \label{prop2}
Let $f$ :\ $U$ $\longrightarrow$ $\mathbb{C}^{n}$, $f(z)=(f_{1}(z),..., f_{n}(z))$, be a finely
holomorphic map on a finely open subset $U\subseteq \mathbb{C}$.
Then the graph $\Gamma _{f}(U)$ of $f$ is a pluripolar subset of
$\mathbb{C}^{n+1}$.
\end{proposition}

\begin{proof} Since the fine topology is
locally connected (Proposition \ref{prop0}), it follows from the
Quasi-Lindel\"{o}f property (Theorem \ref{thm5}) that $U$ has at
most countably many finely connected components. Because a
countable union of pluripolar sets is pluripolar, there is no loss
of generality if we assume that the set $U$ is a fine domain. Let
$a \in U$. According to Theorem \ref{thm8} there exist $V \subset
U$ a finely open fine neighborhood of $a$, and $ \varphi_{j} \in
L^{2}( \mathbb{C}) $, $j=1,\ldots, n$, with compact support such
that $ \varphi_{j} = 0$, $j=1,\ldots, n$, a.e. in $V$ and
$$
f_{j}(z)= \int \limits_{\mathbb{C}} \frac{1}{z-\zeta} \varphi_{j}
(\zeta) d \lambda (\zeta), \ \ z\in V, \ j=1,...,n.
$$
Because of local connectedness, we can assume that $V$ is finely connected. Let $z_{0}
\in V$ and $ 0<\delta <1 $ such that $a \not\in
\overline{D}(z_{0}, \delta)$. Choose a smooth function $ \rho$
such that $ \rho \equiv 1 $ on $D(z_{0}, \delta/2)$ and $\rho
\equiv 0 $ on $ \mathbb{C} \backslash D(z_{0}, \delta)$. Then
$$
\varphi_{j}(\zeta) = \rho(\zeta)\varphi_{j}(\zeta)
+(1-\rho(\zeta))\varphi_{j}(\zeta), \  j=1,...,n.
$$
We set
$$
f_{j}^{1}(z)=\int \limits_{\mathbb{C}} \frac{\rho(\zeta)}{z-\zeta}
\varphi_{j} (\zeta) d \lambda (\zeta), \ \text{and} \ \
f_{j}^{2}(z)=\int \limits_{\mathbb{C}}
\frac{1-\rho(\zeta)}{z-\zeta} \varphi_{j} (\zeta) d \lambda
(\zeta), \  \ j=1,..., n.
$$
It is clear that $f_{j}^{2}$, $j=1,\ldots, n$, is holomorphic on
$D(z_{0}, \delta/2)$ and finely holomorphic on the finely open set
$V \cup D(z_{0}, \delta/2)$. Since usual domains are also finely
connected, $V \cup D(z_{0}, \delta/2)$ is finely connected. Now,
by Lemma \ref{lem2}, the image of $V \cup D(z_{0}, \delta/2)$
under $z \mapsto (z, f_{1}^{2}(z),\ldots,f_{n}^{2}(z))$ is a
pluripolar subset of $\mathbb{C}^{n+1}$. By Josefson's theorem,
there exists a plurisubharmonic function $h
\in\PSH(\mathbb{C}^{n+1})$ ($\not\equiv -\infty$) such that
$$
h(z, f_{1}^{2}(z),\ldots, f_{n}^{2}(z))= - \infty, \ \forall z \in  V
\cup D(z_{0}, \delta/2).
$$
Since $f_{j}^{1}$, $j=1,\ldots,n$, is holomorphic on $\mathbb{C}
\backslash \overline{D}(z_{0}, \delta)$, the function $g$ :\
$\mathbb{C}^{n+1}$ $\longrightarrow$ $\mathbb{C}^{n+1}$, defined
by
$$
g(z, w_{1},\ldots, w_{n})= (z, w_{1}-f_{1}^{1}(z),\ldots, w_{n}-
f_{n}^{1}(z)),
$$
is holomorphic on $\mathbb{C} \backslash \overline{D}(z_{0},
\delta) \times \mathbb{C}^{n}$. Hence $ h\circ g$ is
plurisubharmonic on $\mathbb{C} \backslash \overline{D}(z_{0},
\delta) \times \mathbb{C}^{n}$ and clearly not identically equal
to $- \infty$. Moreover, we have :
$$
h\circ g (z, f_{1}(z),\ldots, f_{n}(z))= h(z,
f_{1}^{2}(z),\ldots,f_{n}^{2}(z))= -\infty , \ \forall z \in V \cap
\mathbb{C} \backslash \overline{D}(z_{0}, \delta).
$$
This proves that the graph $\{(z, f_{1}(z),\ldots,f_{n}(z)) : z \in V
\cap \mathbb{C} \backslash \overline{D}(z_{0}, \delta)\}$ over $V
\cap \mathbb{C} \backslash \overline{D}(z_{0}, \delta)$ is
pluripolar subset of $\mathbb{C}^{n+1}$. Notice that $V \cap
\mathbb{C} \backslash \overline{D}(z_{0}, \delta)$ is a finely
open set containing the point $a$. Again, by Josefson's theorem,
there exists a plurisubharmonic functions $\psi \in$
PSH$(\mathbb{C}^{n+1})$ such that
$$
\psi (z, f_{1}(z),\ldots, f_{n}(z)) = - \infty , \  \forall z \in V
\cap \mathbb{C} \backslash \overline{D}(z_{0}, \delta).
$$
In view of Lemma \ref{lem1} the function $z \longmapsto \psi (z,
f_{1}(z),\ldots, f_{n}(z))$ is finely subharmonic in $U$ or $\equiv -
\infty$. Since it assumes $- \infty$ on the non polar set $V \cap
\mathbb{C} \backslash \overline{D}(z_{0}, \delta)$, it must be
identically equal to $- \infty$ on $U$ by Theorem \ref{thm6}. This
completes the proof of the proposition.
\end{proof}

For convenience of the reader we repeat the statement of our main result, which we will prove subsequently.
\begin{theorem} \label{thm2a} Let $f$ :\ $U$ $\longrightarrow$ $\mathbb{C}^{n}$, $f(z)=(f_{1}(z),\ldots, f_{n}(z))$,
be a finely holomorphic map on a finely open subset $U\subseteq
\mathbb{C}$. Then the image $f(U)$ of $U$ is a pluripolar subset
of $\mathbb{C}^{n}$. Moreover, if $E$ is a non polar subset of
$U$, then the pluripolar hull of $f(E)$ contains $f(U)$.
\end{theorem}

\begin{proof} Without loss of generality we may
assume that $f_{1}$ is not constant and $U$ is a fine domain. It
follows from Theorem \ref{thm10} that one can choose a non empty finely
open subset $W \subseteq U$ of $U$ such that $f_{1}|_{W}$ :\ $W$
$\longrightarrow$ $f_{1}(W)$ is bijective and the inverse function
$f^{-1}_{1}$ is finely holomorphic in the finely open set
$f_{1}(W)$. Now, observe that
$$
f(W)=\{(f_{1}(z),\ldots, f_{n}(z)) : z \in W\} = \{(w,
f_{2}(f_{1}^{-1}(w)),\ldots, f_{n}(f_{1}^{-1}(w))) : w \in
f_{1}(W)\},
$$
where $w = f_{1}(z)$. Since the composition of two finely
holomorphic functions is finely holomorphic (Theorem \ref{thm10}), the map $w \mapsto (f_{2}(f_{1}^{-1}(w)),\ldots,
f_{n}(f_{1}^{-1}(w)))$ is finely holomorphic in $f_{1}(W)$.  Proposition \ref{prop2} applies,
hence the graph
$$\{(w, f_{2}(f_{1}^{-1}(w)),\ldots,
f_{n}(f_{1}^{-1}(w))) : w \in f_{1}(W)\} = f(W)
$$
is a pluripolar subset of $\mathbb{C}^{n}$. Again, Josefson's
theorem ensures the existence of a plurisubharmonic function $h
\in\PSH(\mathbb{C}^{n})$ such that
$$
h(f_{1}(z),\ldots, f_{n}(z)) = -\infty, \forall z \in W.
$$
By Lemma \ref{lem1}, the function $z \mapsto h(f_{1}(z),\ldots,
f_{n}(z))$ is either finely subharmonic or identically equal to
$-\infty$. Since it assumes $-\infty$ on the non polar subset $W
\subset U$, we must have $h(f(U))= -\infty$ by Theorem \ref{thm6}.
Repeating this last argument, the second statement of Theorem
\ref{thm2a} follows. The proof is complete.
\end{proof}

\begin{proof}[Proof of Theorem \ref{thm3}] (1) $\Rightarrow$ (2). According to
Theorem \ref{thm4}, there exists $R > 0$ such that the set $\{z\in
D \backslash A: |f(z)| \geq R \}$ is thin at $z_{0}$. Clearly, the
set $U = \{z\in D \backslash A : |f(z)|< R \} \cup \{z_{0}\}$ is a
finely open neighborhood of $z_{0}$. Since $f$ is bounded in $U
\backslash \{z_{0}\}$ and finely holomorphic in $U \backslash
\{z_{0}\}$, Theorem \ref{thm10} 5) gives that $f$ has a finely
holomorphic extension at $z_{0}$.\\
(2) $\Rightarrow$ (1). Suppose that $f$ has a finely holomorphic
extension $ \tilde{f}$ at $z_{0}$. Clearly, $(D \backslash A) \cup
\{z_{0}\}$ is a finely open neighborhood of $z_{0}$. Since polar
sets do not separate fine domains (Proposition \ref{prop0}) the
set $(D \backslash A) \cup \{z_{0}\}$ is finely connected. Let $h
\in\PSH(D\times\mathbb{C})$ be a plurisubharmonic function such
that $h(z, f(z))=-\infty $, $\forall z \in D \backslash A $.
According to Lemma \ref{lem1}, the function $z \mapsto h(z,
\tilde{f}(z))$ is either finely subharmonic on $(D \backslash
A)\cup \{z_{0}\}$ or $\equiv -\infty$. As it assumes $-\infty$ on
$D \backslash A$, it must be identically equal to $-\infty$ in
view of Theorem \ref{thm6}. Consequently, $(z_{0}, \tilde{f}
(z_{0})) \in (\Gamma_{f})^{\ast}_{D \times \mathbb{C}}$. The last
assertion follows from Theorem 5.10 in \cite{EW04}.
\end{proof}

\section{A Borel-type example}\label{secExample}
We give an example in the spirit of Borel to which the theory
applies. It consists of a finely holomorphic function on a fine
domain, which is a dense subset of $\CC$ with empty Euclidean
interior. Our point is to show that the study of quite natural
series in connection with pluripolarity is fruitfully done in the
framework of fine holomorphy.

\begin{example}

Let $\{a_j\}_{j=1}^\infty$ be a dense sequence in $\CC$ with the property that $|a_j|<j$.
Let $r_j=2^{-j}$.
%\begin{equation}\label{equ1}
%r_j=\min\{\frac{|a_k-a_j|}{2}, 2^{-j}; k=1,\ldots,j-1\}.
%\end{equation}
Then $\cup_{j=1}^\infty B(a_j,r_j)$ has finite area, and its
circular projection $z\mapsto|z|$ has finite length. Next, define
subharmonic functions $g_j(z)=\log|z-a_j|-3j$ and $u_n$ by
\begin{equation}\label{equ2}
u_n(z)=\sum_{j=n}^\infty j^{-3}g_j(z).
\end{equation}
The terms in the sum of \eqref{equ2} are subharmonic and they are
negative for $|z|<k$ as soon as $j>k$. Hence $u_n$ represents a
subharmonic function. Let
$D=\left(\cup_{n}\{u_n>-10\}\right)\setminus\{a_1,a_2,\ldots\}$.
We claim that $D = \{u_1>-\infty\}$. Indeed, let $z_0 \in \cup_{n}
\{u_n
> -\infty \} \setminus\{a_1,a_2,\ldots\}$. Then there exists a natural number $k$ such
that $|z_0| < k $ and $u_k > -\infty $. Since, as mentioned
before, all the terms of the series $u_k (z_0)$ are negative, a
suitable tail, say $u_N (z_0)$, will be very close to $0$. In
other words, $z_0 \in \{u_N > -10\}$. Hence $z_0 \in D$ and
consequently $D= \cup_{n} \{u_n
> -\infty \} \setminus\{a_1,a_2,\ldots\}$. Therefore,
$$
\mathbb{C} \backslash D = \cap_{n=1}^{\infty} \{u_k =- \infty \}
\cup \{a_1,a_2,\ldots\}.
$$
Since $\{u_{k_1} =- \infty \} \backslash
\{a_1,a_2,\ldots\}=\{u_{k_2} =- \infty \} \backslash
\{a_1,a_2,\ldots\}$ for any natural numbers $k_1$ and $k_2$, we
conclude that
$$
\mathbb{C} \backslash D = \{u_1 =- \infty \} \cup
\{a_1,a_2,\ldots\}=\{u_1 =- \infty \}.
$$
This proves the claim. In particular, $D$ is, by Proposition
\ref{prop0}, a fine domain.

For every $j$ there exists $0<c_j<1$ such that if $|z-a_j|<c_j$,
then for $n\le j$, $u_n(z)<-11$. Indeed, $\sum_{k>j}
k^{-3}g_k(z)<0$, while
$$\sum_{k=n}^{j-1} k^{-3}g_k(z)<\log j\sum_{k=n}^{j-1}k^{-3}<10\log j.$$
So it suffices to take $c_j=j^{-11j^3}$.

Next we define a function on $D$ by
\begin{equation}\label{equ3}
f(z)=\sum_{j=1}^\infty\frac{c_j}{2^j(z-a_j)},
\end{equation}
We claim that the function $f$ is finely holomorphic on $D$.
Indeed, let $z_0\in D$. For every $m$  a suitable tail of the
series of $f$ in \eqref{equ3} is uniformly convergent on the
compact set  $K=\{|z|\le 2|z_0|\}\setminus \cup_{j\ge
m}B(a_j,c_j)$. Now if $z_0\in D$, then $z_0$ belongs to the finely
open set $\{u_m>-10\}$ for some $m$. Hence, for all $j\ge m$ we
have $|z_0-a_j|>c_j$, and $K$ is a fine neighborhood of $z_0$.
\end{example}

Application of Theorem \ref{thm2} gives us that the graph of
$\Gamma_f(D)$ of $f$ over $D$ is a pluripolar set. The theorem
also shows that for  a set of positive capacity $E\subset D$,
e.g., a circle in $\{u_1>-10\}$,
$$\Gamma_f(D)\subset \left(\Gamma_f(E)\right)_{\CC^2}^*.$$

Even for this example there are many questions left. We have no
description of the maximal domain $D_0$  to which $f$ extends as a
finely holomorphic function, and we don't know if $\Gamma_f(D_0)=
\left(\Gamma_f(E)\right)_{\CC2}^*$, as one may expect in view of
\cite{EW04}.
%Note that one can argue as in the proof of
%(2)$\Rightarrow$ (1) of Theorem 1.3 to show that
%$\left(\Gamma_f(E)\right)_{\CC^2}^*=\Gamma_f(D_0).$

%\begin{center}
\section{Pluripolarity of finely analytic varieties}
%\end{center}

In this section we will extend Theorem \ref{thm1} to the case
where $f$ is a function of several complex variables. To do this,
we will first define finely plurisubharmonic functions and finely
holomorphic functions of several variables.
\begin{definition}
1) The {\em pluri-fine topology} on a domain $\Omega$ is the weakest topology
that makes all plurisubharmonic functions continuous.\\
2) A function $f$ on a pluri-fine domain $\Omega$ is called {\em
finely plurisubharmonic} if it is upper semi-continuous (in the
pluri-fine topology) and if the restriction of $f$ to any complex
line $L$ is
finely subharmonic or identically $-\infty$ on any finely connected component of $\Omega\cap L$.\\
3) A function $f$ on a pluri-fine domain is called {\em finely
holomorphic} if every point $z\in\Omega$ has a Euclidean compact,
pluri-fine neighborhood $K$, such that $f\in H(K)$.
\end{definition}
Here $H(K)$ stands for the uniform closure on $K$ of the algebra of holomorphic functions in a neighborhood of $K$.

Proposition \ref{prop0} remains valid in the pluri-fine setting.
Items 1) and 2) of it were proved in \cite{E-W1} and 3) in
\cite{E-W2}. It would be interesting to know if in general finely
pluripolar sets are pluripolar, as is the case in dimension 1.
However, for our purposes here, a weaker result that we proved in
\cite{E-W2}, suffices.

\begin{proposition} \label{prop5}Let $E$ be a pluri-finely open subset of a pluri-fine domain $D \in \CC^n$.
Suppose that $h$ is finely plurisubharmonic on $D$ and
$h|_E=-\infty$. Then $h\equiv-\infty$ on D. In other words,
pluri-finely open sets are not finely pluripolar.
\end{proposition}

\begin{lemma}\label{flem} Let $U\subseteq\CC^{n}$ be a pluri-fine domain, and let $f$ :\
$U$ $\longrightarrow$ $\CC$ be a finely holomorphic function.
Suppose that $h$ :\ $\CC^2$ $\longrightarrow$ $[-\infty,\ +\infty
[$ is a plurisubharmonic function. Then the function $z \mapsto
h(z,f(z))$ is finely plurisubharmonic on $U$.
\end{lemma}

The proof is exactly the same as the proof of Lemma \ref{lem1}.

\begin{theorem}\label{thm11} Let $f$ be a finely holomorphic function in a pluri-fine domain $U
\subset \CC^{n}$. Suppose that for some pluri-finely open subset
$V\subset U $ the graph $ \Gamma_{f}(V)$ of $f$ over $V$ is
pluripolar in $\CC^{n+1}$. Then the graph $ \Gamma_{f}(U)$ of $f$
is pluripolar in $\CC^{n+1}$. Moreover, $ \Gamma_{f}(U)\subset
(\Gamma_{f}(V))^{\ast}_{\CC^{n+1}}$.
\end{theorem}
\begin{proof} By Josefson's theorem there
exists $h\in \PSH(\CC^{n+1})$ ($\not\equiv -\infty$) such that
$h(z, f(z))= -\infty$, $\forall z\in V$. In view of Lemma
\ref{flem} and Proposition \ref{prop5}, the function $h(z, f(z))$
is identically $-\infty$ in $U$. It follows at once that $
\Gamma_{f}(U)$ is pluripolar and contained in the pluripolar hull
of $\Gamma_{f}(V)$.
\end{proof}

As a corollary we obtain a generalization of Theorem \ref{thm1}. We keep the notation of Theorem \ref{thm11}.
\begin{corollary} \label{cor1} Suppose that $U$ contains a Euclidean ball $B$.
 Then $ \Gamma_{f}(U)$ is pluripolar and $ \Gamma_{f}(U)\subset
(\Gamma_{f}(B)))^{\ast}_{\CC^{n+1}}$.
\end{corollary}
\begin{proof} It follows from Theorem \ref{thm9} that $f$ is holomorphic on the intersection
of $B$ with any complex line. Hence $f$ is holomorphic on $B$ and
$\Gamma_{f}(B)$ is pluripolar. Now Theorem \ref{thm11} applies.
\end{proof}

\begin{remark}
 Theorem \ref{thm11} and Corollary \ref{cor1} only explain for a small part the propagation
of pluripolar hulls. E.g., in the case of Corollary \ref{cor1}
take $B$ the unit ball and consider the function
$g(z)=f(z)(z_1-z_2^2)$. Then, whatever the extendibility
properties of $f$ may be, the pluripolar hull of the graph of $g$
will contain the set $\{z_1=z_2^2\}$.
\end{remark}

\section{Concluding remarks and open questions}

%Let $U\subset
%\mathbb{C}$ be a bounded fine domain and let $f$ :\ $U$
%$\longrightarrow$ $\mathbb{C}^{n}$ be a finely holomorphic map. We
%call $f(U)$ a \emph{finely holomorphic curve}.
We now discuss some applications and open problems. Let $E\subset
\mathbb{C}^{n}$ be a pluripolar set and $E^{*}_{\mathbb{C}^{n}}$
its pluripolar hull. It follows from the arguments used before
that if $E$ hits a finely holomorphic curve $f(U)$ in some non
"small" set, then $E^{*}_{\mathbb{C}^{n}}$ contains all the points
of $f(U)$. Namely, we have the following.
\begin{proposition} \label{prop3}
Let $f$ :\ $U$ $\longrightarrow$ $\mathbb{C}^{n}$ be a finely
holomorphic map on a bounded fine domain $U \subset \mathbb{C}$
and let $E \subset \mathbb{C}^{n}$ be a pluripolar set. If
$f(U)\cap E \neq \emptyset$ and $f^{-1}(f(U)\cap E )$ is
non-polar, then $f(U)\subset E^{*}_{\mathbb{C}^{n}}$.
\end{proposition}

\begin{proof} Let  $h \in$ PSH$(\mathbb{C}^{n})$ be a
plurisubharmonic function such that $h(z)= - \infty$, $\forall z
\in E$. By Lemma \ref{lem1}, $h\circ f$ is either finely
subharmonic on $U$ or $\equiv -\infty$. As it assumes $-\infty$ on
$f^{-1}(f(U)\cap E)$, it must be, by Theorem \ref{thm6},
identically $-\infty$ on $U$. We have therefore $f(U)\subset
E^{*}_{\mathbb{C}^{n}}$.
\end{proof}

The conclusion of the above proposition remains valid if one
assumes that $E$ contains merely the "boundary of a finely
holomorphic curve".
\begin{proposition} \label{prop4} Let $f$ and $E$ be as above. If $f$ extends by fine continuity to the fine
boundary $\partial _{f}U$ of $U$ and $f(\partial_{f}U) \subset E$,
then $f(U)\subset E^{*}_{\mathbb{C}^{n}}$.
\end{proposition}

\begin{proof} Let  $h \in$ PSH$(\mathbb{C}^{n})$ be
plurisubharmonic function such that $h(z)= - \infty$, $\forall z
\in E$. Let $a \in \partial _{f}U$. By assumption, $f$ has a fine
limit at $a$. Using Cartan's theorem (cf. \cite{H69}, Theorem
10.15), one can easily find a finely open neighborhood $V_{a}$ of
$a$ such that the usual limit $\lim_{z \to a,z\in V_{a}\cap
U}f(z)$ exists and is equal to $f(a)$. Let $M > 0$. Since $h$ is
upper semicontinuous, the set $\{ z \in \mathbb{C}^{n}: h(z)< -
M\}$ is open. As $ f(a) \in \{ z \in \mathbb{C}^{n}: h(z)< - M\}$,
one can find a positive number $\delta_{a} > 0$ such that
$$
f(w) \in \{ z \in \mathbb{C}^{n}: h(z)< - M\}, \ \forall w \in
B(a, \delta_{a}) \cap V_{a} \cap U,
$$
where $B(a, \delta_{a})$ is the disk with center $a$ and radius
$\delta_{a}$. Consequently
$$
\f-limsup_{z \to a,z\in U}h(f(z)) \leq \f-limsup_{z \to a,z\in
V_{a} \cap U}h(f(z)) < -M, \ \forall a \in \partial _{f}U
$$
where $\f-limsup$ denotes the limit with respect to the fine
topology. As $h\circ f$ is a finely hypoharmonic function on $U$
(see Lemma \ref{lem1} and its proof), the fine boundary maximum principle
(cf. \cite{Fu74}, Theorem 2.3) shows that $h\circ f(z) < -M$, $
\forall z \in U$. Since $M$ was arbitrary, we conclude that
$h\circ f(U)= - \infty$. This proves the proposition.
\end{proof}

Our results reveal a very close relationship between the
pluripolar hull of the graph of a holomorphic function and the
theory of finely holomorphic functions (see also \cite{EJ}). This
leads
naturally to the following fundamental problem.  \\
\textbf{Problem 1}. Let $f$ :\ $\Omega $ $\longrightarrow$
$\mathbb{C}$ be a holomorphic function on a simply connected open
subset $\Omega \subset \mathbb{C}$. Suppose that the graph
$\Gamma_{f}(\Omega)$ of $f$ over $\Omega $ is not complete
pluripolar. Must then
$(\Gamma_{f}(\Omega))^{*}_{\mathbb{C}^{2}}\backslash
\Gamma_{f}(\Omega)$ have a fine analytic structure? i.e., Let
$z\in (\Gamma_{f}(\Omega))^{*}_{\mathbb{C}^{2}}\backslash
\Gamma_{f}(\Omega)$. Must there exist a finely holomorphic curve
passing through $z$ and contained in
$(\Gamma_{f}(\Omega))^{*}_{\mathbb{C}^{2}}\backslash
\Gamma_{f}(\Omega)$?\\
Obviously, a positive answer to the above problem would, in
particular, solve the following problem posed
in \cite{EJ}.\\
\textbf{Problem 2}. Let $f$ be a holomorphic function in the unit
disc $\mathbb{D}$. Suppose that
$(\Gamma_{f}(\mathbb{D}))^{*}_{\mathbb{C}^{2}}$ is the graph of
some function $\mathcal{F}$. Is  $\mathcal{F}$ a finely
holomorphic continuation of $f$?

It was proved in \cite{CLP} that one can not detect
"pluripolarity" via intersection with one dimensional complex
analytic varieties. Since there are, roughly speaking, much more
finely holomorphic curves in $\mathbb{C}^{n}$
than analytic varieties, one can naturally pose the following\\
\textbf{Problem 3}. Let $K$ be a compact set in $\mathbb{C}^{n}$
and suppose that $f^{-1}(K\cap f(U))$ is a polar subset of $U$ (or
empty) for any finely holomorphic curve $f$ :\ $U$
$\longrightarrow$ $\mathbb{C}^{n}$. Must $K$ be a pluripolar
subset of $\mathbb{C}^{n}$?

\noindent Kd\textsc{V Institute for Mathematics, University of
Amsterdam, Plantage Muidergracht, 24, 1018 TV, Amsterdam,
The Netherlands}\\
janwieg@science.uva.nl \\  smarzgui@science.uva.nl

\noindent\textsc{Institute of Mathematics, Jagiellonian
University, Reymonta 4, 30-059 Krak\'{o}w,
Poland}\\
Armen.Edigarian@im.uj.edu.pl


\begin{thebibliography}{99}

%\bibitem{DG74} Debiard, A., Gaveau, B.: Potential fin et
%alg\`ebres de fonctions analytiques, I, {\em J. Functional Anal.}
%{\bf 16} (1974), 289--304, II, {\bf 17} (1974), 296--310.
\bibitem{CLP} Coman. D., Levenberg, N., Poletsky, E. A.: Smooth submanifolds intersecting any analytic curve
in a discrete set, {\em Math. Ann.}, {\bf 332} (2005), 55--65.

\bibitem{DO84} Doob, J. L.: "Classical Potential Theory and Its Probabilistic Counterpart,
" Springer-Verlag, Berlin, 1984.

\bibitem{EW03}Edigarian, A., Wiegerinck, J.: The pluripolar hull of the graph of a holomorphic
function with polar singularities, {\em Indiana Univ. Math. J.},
{\bf 52} (2003), 1663--1680.

\bibitem{EW04}Edigarian, A., Wiegerinck, J.: Determination of the pluripolar hull of graphs
of certain holomorphic functions, {\em Ann. Inst. Fourier.
Grenoble}, {\bf 54} (2004), no. 6, 2085--2104.

\bibitem{E} Edlund, T.: Pluripolar sets and pluripolar hulls, {\em Uppsala Dissertations in Mathematics} {\bf
41} (2005).

\bibitem{EJ} Edlund, T., J\"oricke, B.: The pluripolar hull of a graph and fine analytic
continuation, {\em Ark. Math.} {\bf 44 } (2006), no. 1, 39--60.

\bibitem{E-W1} El Marzguioui, S. Wiegerinck, J.: The pluri-fine topology is locally connected.
{\em Potential Anal} {\bf 25} (2006), no. 3, 283--288.

\bibitem{E-W2} El Marzguioui, S., Wiegerinck, J.: Connectedness in the Pluri-fine Topology, preprint 2008.

%\bibitem{Fu71} Fuglede, B.: Connexion en topologie fine et balayage des mesures, {\em
%Ann. Inst. Fourier Grenoble} {\bf 21.3} (1971), 227--244.

\bibitem{Fu72} Fuglede, B.: Finely harmonic functions, {\em Springer Lecture Notes in Mathematics,} {\bf 289},
Berlin-Heidelberg-New York, 1972.

\bibitem{Fu74} Fuglede, B.: Fonctions harmoniques et fonctions finement harmoniques, {\em
Ann. Inst. Fourier Grenoble.} {\bf 24.4} (1974), 77--91.

\bibitem{Fu75} Fuglede, B.: Asymptotic paths for subharmonic functions, {\em Math. Ann,} {\bf 213}, (1975), 261-274.

\bibitem{Fu76} Fuglede, B.: Finely harmonic mappings and finely holomorphic functions, {\em
Ann. Acad. Sci. Fennic\ae,} {\bf 2} (1976), 113--127.

\bibitem{Fu81} Fuglede, B.: Sur les fonctions finement holomorphes, {\em
Ann. Inst. Fourier.} {\bf 31.4} (1981), 57--88.

%\bibitem{Fu80} Fuglede, B.: Fine topology and finely holomorphic functions, {\em
%"Proc. 18th Scand. Congr. Math. Aarhus 1980",}
%Birkh\"auser, (1981), pp. 22--38.

%\bibitem{Fu82} Fuglede, B.: Localisation in fine potential theory and uniform approximation by subharmonic
%functions, {\em
%J. Functional Anal.} {\bf 49} (1982), 57--72.

\bibitem{Fu88} Fuglede, B.: Finely holomorphic functions- a survey, {\em
Revue Roumaine Math. Pures Appl.} {\bf 33} (1988), 283--295.

\bibitem{H69} Helms, L. L.: Introduction to potential theory, {\em Pure and Applied Mathematics}, Vol
XXII. Wiley-Interscience, New York, 1969.

\bibitem{Jo78} Josefson, B.: On the equivalence between locally polar and globally polar sets for
plurisubharmonic functions on $(\mathbb{C}^{n})$ ,  {\em
Ark. Math.} {\bf 16} (1978), 109--115.

\bibitem{K91} Klimek, M.: Pluripotential Theory, London Mathematical Society Monographs, 6,
Clarendon Press, Oxford, 1991.

\bibitem{Le88} Levenberg, N.: On an example of Wermer, {\em
Ark. Math.} {\bf 26 } (1988), no. 1,  155--163.

\bibitem{LP99} Levenberg, N., Poletsky, E. A.: Pluripolar hulls, {\em
Mich. Math. J.} {\bf 46 } (1999), 151--162.

\bibitem{LS} Levenberg, N., Slodkowski, Z.: Pseudoconcave pluripolar sets in $\mathbb{C}^{2}$ , {\em
Math. Ann.} {\bf 312 } (1998), 429--443.

%\bibitem{Ly80} Lyons, T. J.: Finely holomorphic functions, {\em
%J. Functional Anal.} {\bf 37}  (1980), 1--18.

%\bibitem{N78}Nguyen-Xuan-Loc.: Sur la th\'{e}orie des fonctions finement holomorphes. {\em
%Bull. Sci. Math.}, {\bf 102} (1978), 271--308.

\bibitem{Sh03} Shcherbina, N. V.: Pluripolar graphs are
holomorphic, {\em Acta. Math.}, {\bf 194} (2005), 203--216.

\bibitem{Si03} Siciak, J.: Pluripolar sets and pseudocontinuation, in: Complex Analysis
and Dynamical Systems II (Nahariya, 2003), Contemp. Math. 382,
{\em,Amer. Math. Soc}., Providence, RI, 2005, 385--394.  {\bf 22}
(2005),  195--206.

\bibitem{We82} Wermer, J.: Polynomially convex hulls and analyticity , {\em
Ark. Math.} {\bf 20 } (1982), 129--135.

\bibitem{Wi00} Wiegerinck, J.: Graphs of holomorphic functions with isolated singularities are complete
pluripolar, {\em Michigan Math. J.}, {\bf 47} (2000), 191--197.

\bibitem{Ze89} Zeriahi, A.: Ensembles pluripolaires exceptionnels pour la croissance partielle des fonctions
holomorphes, {\em Ann. Polon. Math.} {\bf 50 } (1989), 81--91.

\bibitem{Zw05} Zwonek, W.: A note on pluripolar hulls of graphs of Blaschke products, {\em
Potential Analysis,} {\bf 22} (2005), 195--206.

\end{thebibliography}
\end{document}